\documentclass[12pt,eqno]{article}
\usepackage{amsmath,amsthm,amscd,amssymb}
\usepackage{latexsym}
\newtheorem{remark}{Remark}[section]
\newcommand{\bremark}{\begin{remark} \em}
\newcommand{\eremark}{\end{remark} }

\begin{document}
\parindent 12pt
\renewcommand{\theequation}{\thesection.\arabic{equation}}
\renewcommand{\baselinestretch}{1.15}
\renewcommand{\arraystretch}{1.1}
\def\disp{\displaystyle}
\title{\large Positve solutions of inhomogeneous Kirchhoff type\\ equations with indefinite data\\
\footnotetext{E-mail addresses: qiuyidai@aliyun.com, 2432295269@qq.com.\\
\indent\indent This work is supported by NNSFC (Grant:No.11671128).}}
\author{{\small Aolin Chen,\      Qiuyi Dai}\\
{\small\emph{ LCSM, College of Mathematics and Statistics, Hunan Normal University,}}\\
\small\emph{{Changsha Hunan 410081,P. R. China}}\\
\date{}}
\maketitle

\abstract
{\small Let $\Omega$ be a bounded domain in $\mathbb{R}^N$ with smooth boundary $\partial\Omega$. Denote by $\mathcal{M}$ the subset of $C^1(\overline{\Omega})\backslash \{0\}$ such that for any $f(x)\in\mathcal{M}$ the following problem
\begin{equation}\label{eq1}\begin{split}
\left\{\begin{array}{ll}
-\Delta u=f(x) & x\in \Omega, \\
u\geq 0 & x\in\Omega, \\
u=0 & x\in \partial\Omega,
\end{array}
\right.
\end{split}\end{equation}
has a solution.
Assume that $b>0,p>1$ and $\lambda>0$. We consider Dirichlet problem of inhomogeneous Kirchhoff type equation
\begin{equation}\label{eq2}\begin{split}
\left\{\begin{array}{ll}
-(1+b\|\nabla u\|^{2\alpha}_2)\Delta u=u^{p}+\lambda f(x) & x\in\Omega, \\
u> 0 & x\in\Omega, \\
u=0 & x\in \partial\Omega,
\end{array}
\right.
\end{split}\end{equation}
where $\alpha \in {(0,{\frac{2^*-1}{2}})}$ with $2^*=+\infty$ for $N=2$, and $2^*={\frac{N+2}{N-2}}$ for $N\geq 3$. 

Main results we proved in the present paper can be summarized as

\noindent {\bf (i)}\ \ If $1<p<2\alpha+1$, then, for any $\lambda>0$ and $f(x)\in\mathcal{M}$, problem (\ref{eq2}) has at least one solution.

\noindent {\bf (ii)}\ \ If $1<p<2\alpha+1$ and $b>b_0$ for some positive number $b_0$ given by (\ref{eq7}) in Section 1, then problem (\ref{eq2}) is solvable if and only if $f(x)\in\mathcal{M}$. Moreover, the solution is unique for $\lambda$ small enough.

\noindent {\bf (iii)}\ \ If $2\alpha+1<p<2^*$ and $f(x)\in\mathcal{M}$, then problem (\ref{eq2}) has at least two solutions for $\lambda$ small enough and has no solution for $\lambda$ large enough.

\noindent {\bf (iv)}\ \ If $p>2^*$, then problem (\ref{eq2}) has at least one solution for $\lambda$ small enough if and only if $f(x)\in\mathcal{M}$, and has no solution for $\lambda$ large enough.

Compared to the semilinear case (that is the case $b=0$), the appearance of the nonlocal term $b\|\nabla u\|^2$ in Kirchhoff type equations changes tremendously the profile of the solution set in the case $1<p<2\alpha+1$. For more detailed explanation, see Remark 1.4 in Section 1.}


\section*{1. Introduction}

\setcounter{section}{1}

\setcounter{equation}{0}

\noindent
Let $\Omega$ be a bounded domain in $\mathbb{R}^N$ with smooth boundary $\partial\Omega$, and $ f(x)\in C^1(\overline{\Omega})\backslash \{0\}$. For any $1\leq q\leq\infty$, we use $L^q(\Omega)$ to denote the standard Lebesgue's space endowed with norm $\|\centerdot\|_q$. In this paper, we consider the following Dirichlet problem of inhomogeneous Kirchhoff type equation
\begin{equation}\label{eq3}\begin{split}
\left\{\begin{array}{ll}
-(1+b\|\nabla u\|^{2\alpha}_2)\Delta u=|u|^{p-1}u+\lambda f(x) & x\in\Omega, \\
u=0 & x\in\partial\Omega,
\end{array}
\right.
\end{split}\end{equation}
where $b>0, p>1$, $\lambda >0$, and $0<\alpha <{\frac{2^*-1}{2}}$ with $2^*=+\infty$ for $N=2$, and $2^*={\frac{N+2}{N-2}}$ for $N\geq 3$. 

Since the differential equation in problem (\ref{eq3}) contains an integral over $\Omega$, it is no longer a pointwise identity. Therefore, it is often called nonlocal problem. Nonlocal boundary value problems like problem (\ref{eq3}) model several physical and biological systems where $u$ describes a process which depend on the average of itself, such as the population density. We refer the reader to \cite{Va, AC, AM, CL, CR} for some related works. Concerning problem (\ref{eq3}) itself, the prototype of it is the Kirchhoff wave equation
which was proposed by Kirchhoff in \cite{Kir} as an extension of the classical D'Alembert's wave equation, by considering the effect of the changing in the length of the string during the vibration. For more mathematical and physical background of Kirchhoff equations, we refer to \cite{Aro, CDCS, DS, HZ} and the references cited there in. 

In the case $b=0$, problem (\ref{eq3}) is reduced to the following well studied semilinear problem
\begin{equation}\label{eq4}\begin{split}
\left\{\begin{array}{ll}
-\Delta u=|u|^{p-1}u+\lambda f(x) & x\in\Omega, \\
u=0 & x\in\partial\Omega.
\end{array}
\right.
\end{split}
\end{equation}

To our best knowledge, the study of problem (\ref{eq4}) was initiated by \cite{AB} in which A. Bahri and H. Berestycki tried to find infinitely many nontrivial solutions by perturbation method. Since then, problem  (\ref{eq4}) has attracted many attentions, see for example \cite{AL, Tan, Ra} etc.  What we emphases here are positive solutions of problem (\ref{eq4}). In this respect, many authors have made their contributions under the assumption that $f(x)\geq 0$, see for example \cite{L, D, DL}. The condition $f(x)\geq 0$ has been improved in \cite{DG, DY, DP} by Q. Y. Dai, Y. G. Gu, J. F. Yang and L. H. Peng. To recall the results of  \cite{DG, DY, DP}, We denote by $\mathcal{M}$ the subset of $C^1(\overline{\Omega})\backslash \{0\}$ such that for any $f(x)\in\mathcal{M}$ the following problem
\begin{equation}\label{eq5}\begin{split}
\left\{\begin{array}{ll}
-\Delta u=f(x) & x\in \Omega, \\
u\geq 0 & x\in\Omega, \\
u=0 & x\in \partial\Omega,
\end{array}
\right.
\end{split}
\end{equation}
has a solution. Obviously, $\mathcal{M}$ includes sign-changing function. With the notation $\mathcal{M}$, main results of \cite{DG, DY, DP} can be summarized as
\vskip 0.1in
{\em {\bf Theorem I} Assume that $f(x)\in C^1(\overline{\Omega})\backslash \{0\}$. Then the following statements hold.

\noindent {\bf (i)}\ \ If $1<p<2^*$ and $f(x)\in\mathcal{M}$, then there exists a positive number $\lambda_{f}<+\infty$ such that  problem (\ref{eq4})  has at least two positive solutions for any $\lambda\in(0,\lambda_{f})$, and has no positive solution for $\lambda>\lambda_f$.

\noindent {\bf (ii)}\ \ If $p>2^*$ and $\Omega$ is starshaped, then there exists a positive number $\lambda_{f}<+\infty$ such that  problem (\ref{eq4})  has at least one positive solution for any $\lambda\in(0,\lambda_{f})$ if and only if $f(x)\in\mathcal{M}$, and has no positive solution for $\lambda>\lambda_f$.
}

It is worth pointing out here that sub-supersolution method plays an important role in the study of semilinear problem.

Back to the Kirchhoff type equations (that is the case $b>0$), it attracts more and more attentions in the recent years. See for example \cite{ Azz, AZ, C, CKW,  Na1,Na2, GF, ST, PZ, XC, ZT, CDCS, Ch, DS, ZSN, DLS, Sh, HZ, WHL, LLS, ZS}.  Most literatures available so far are concerning with ground state solutions for homogenous Kirchhoff equations.
However, it is worth mentioning that N. Azzouz and A. Bensedik [31] have studied in \cite{ABe} the following inhomogenous problem
\begin{equation}\label{eq6}\begin{split}
\left\{\begin{array}{ll}
-M(\|\nabla u\|^{2}_2)\Delta u=|u|^{p-1}u+\lambda f(x)& \mbox{in}\ \ \Omega, \\
u=0 & \mbox{on}\ \ \partial\Omega,
\end{array}
\right.
\end{split}
\end{equation}
where $p\in(0,1)\cup(1,2^*), \lambda>0$. 

By making use of sub-supersolution method, they proved that if $M(s)$ satisfies the following conditions:\\
\indent $(M_0)$ $M(s)$ is a continuous and for any $s>0$, $M(s)\geq m_0$ for some $m_0>0$,\\
\indent $(M_1)$ $M(s)$ is a nonincreasing function,\\
\indent $(M_2)$ The function $H(s)=sM(s^2)$ is increasing,\\
then, for any $f(x)\in\mathcal{M}$, there are positive numbers $\lambda^1_f, \lambda_f^2<+\infty$ such that problem (\ref{eq6}) has at least one nonnegative solution for $0<\lambda<\lambda^1_f$, and has no nonnegative solution for $\lambda>\lambda^2_f$

Using the notation $M(s)$ of N. Azzouz and A. Bensedik, we have $M(s)=1+bs^{2\alpha}$ in our problem (\ref{eq3}). This obviously beyond the consideration of \cite{ABe}. Moreover, since $M(s)$ is increasing and unbounded in our problem, the comparison principle may cease to validate (see \cite{GFA}), and sub-supersolution method is no longer available for Kirchhoff type equation itself. Therefore, some new ideas are needed for finding positive solutions of problem (\ref{eq3}) when the data changes sign and $p$ is supercritical. Next, we are going to state our main results of the present paper. To this end, we fix some notations first.

Let $H^1_0(\Omega)$ be the standard Sobolev space and $S(\Omega)$ be the Sobolev constant defined by
$$S(\Omega)=\inf\limits_{u\in H^1_0(\Omega)\backslash\{0\}}\frac{\|\nabla u\|^2_2}{\|u\|^2_{p+1}}.$$
Set $\gamma=2\alpha+1-p$ and $l=S^{\frac{p+1}{2}}(\Omega)$. For $1<p<2\alpha+1$, we introduce a positive constant $b_0$ by the following formula:

\begin{equation}\label{eq7}\begin{split}
\begin{array}{ll}
b_0=(p-1)\gamma^{\frac{\gamma}{p-1}}(2\alpha l)^{-\frac{2\alpha}{p-1}}.
\end{array}
\end{split}
\end{equation}
Bearing above notations in mind, we can express our main results of this paper in the following theorems.
\vskip 0.1in
{\em {\bf Theorem 1.1}\ \ If $1<p<2\alpha+1$ and $f(x)\in\mathcal{M}$, then problem (\ref{eq3}) has at least one positive solution for any $\lambda>0$.}
\vskip 0.1in
{\em {\bf Theorem 1.2}\ \ If $1<p<2\alpha+1$ and $b>b_0$, then problem (\ref{eq3}) has positive solution for any $\lambda>0$ if and only if $f(x)\in\mathcal{M}$. Moreover, the solution is unique for $\lambda$ small enough if in addition $\alpha\geq\frac{1}{2}$.} 
\vskip 0.1in
{\em {\bf Theorem 1.3}\ \ If $2\alpha+1<p<2^*$ and $f(x)\in\mathcal{M}$, then there are two positive constants $\lambda_f, \Lambda_f<+\infty$ such that problem (\ref{eq3}) has at least two positive solutions for $\lambda\in (0, \lambda_f)$, and has no nonnegative solution for $\lambda>\Lambda_f$.}  
\vskip 0.1in
{\em {\bf Remark 1.4}\ \ From Theorem I (i), Theorem 1.1 and Theorem 1.2, we see that the appearance of the nonlocal term $b\|\nabla u\|_2^{2\alpha}$ in Kirchhoff type equation changes the profile of solution set in two aspects when $1<p<2\alpha+1$. One is that the positive solvability of semilinear problem needs a finite restriction on the parameter $\lambda$,whereas Kirchhoff type equation is always positively solvable for any positive parameter $\lambda$; the other one is that semilinear problem has always two positive solutions for small parameter $\lambda$, whereas Kirchhoff type equation has only one positive solution for small parameter $\lambda$ and large $b$ when $\alpha\geq\frac{1}{2}$.} 
\vskip 0.1in
{\em {\bf Theorem 1.5}\ \ If $p>2^*$ and $\Omega$ is starshaped, then there are two positive constants $\lambda_f, \Lambda_f<+\infty$ such that problem (\ref{eq3}) has positive solution for any $\lambda\in(0, \lambda_f)$ if and only if $f(x)\in\mathcal{M}$, and has no positive solution for $\lambda>\Lambda_f$.} 
\vskip 0.1in 
{\em {\bf Remark 1.6}\ \ If not specially declared, all solutions of this paper are in classical sense.}
\vskip 0.1in
The rest of the paper is organized as follows. The case $1<p<2\alpha+1$ is discussed in Section 2. The discussion of the case $2\alpha+1<p<2^*$ is placed in Section 3. The last Section 4 devotes to discuss the case $p>2^*$.

\section*{2. The case $1<p<2\alpha+1$}

\setcounter{section}{2}

\setcounter{equation}{0}

\noindent

Keeping notations $\mathcal{M}$, $\gamma$, and $b_0$ of the previous section in use, we study the case $1<p<2\alpha+1$ in this section. The main results we will prove are  following Theorems.
\vskip 0.1in
{\em {\bf Theorem 2.1}\ \ If $1<p<2\alpha+1$ and $f(x)\in\mathcal{M}$, then problem (\ref{eq3}) has at least one positive solution for any $\lambda>0$.}
\vskip 0.1in
{\em {\bf Theorem 2.2}\ \ If $1<p<2\alpha+1$, and $b>b_0$, then problem (\ref{eq3}) has positive solution for any $\lambda>0$ if and only if $f(x)\in\mathcal{M}$. Moreover, the solution is unique for $\lambda$ small enough if in addition $\alpha\geq\frac{1}{2}$.} 
\vskip 0.1in

To prove Theorem 2.1, we need a result about the solvability of the following problem

\begin{equation}\label{eq21}\begin{split}
\left\{\begin{array}{ll}
-(1+b\|\nabla u\|^{2\alpha}_2)\Delta u=\lambda f(x) &  x\in\Omega, \\
u\geq 0 & x\in\Omega,\\
u=0 &  x\in\partial\Omega,
\end{array}
\right.
\end{split}
\end{equation}
where $b>0,\alpha>0,\lambda>0$. Which can be stated as
\vskip 0.1in

{\em {\bf Lemma 2.3}\ \ Problem (\ref{eq21}) is solvable if and only if $f(x)\in\mathcal{M}$.} 

{\bf Proof:}\ \ On one hand, if $u$ is a solution of problem (\ref{eq21}), then it is easy to check that $v=\frac{1+b\|\nabla u\|^{2\alpha}_2}{\lambda}u$ is a solution of the following problem
\begin{equation}\label{eq22}\begin{split}
\left\{\begin{array}{ll}
-\Delta v=f(x) & x\in\Omega, \\
v\geq 0 & x\in\Omega,\\
v=0 & x\in\partial\Omega.
\end{array}
\right.
\end{split}\end{equation}
Hence, $f(x)\in\mathcal{M}$.

On the other hand, if $f(x)\in\mathcal{M}$, then problem (\ref{eq22}) has a solution $v(x)$. Based on the observation of the above paragraph, we can find a solution of problem (\ref{eq21}) with the form $u_\beta=\frac{\lambda}{1+b\beta^{\alpha}}v$. It is easy to check that $u_\beta$ is indeed a solution of problem (\ref{eq21}) provided that $\beta$ is a positive solution of the following algebraic equation
$$by^{\alpha+\frac{1}{2}}+y^{\frac{1}{2}}-\lambda \|\nabla v\|_2=0. $$
Noting that $h(y)=by^{\alpha+\frac{1}{2}}+y^{\frac{1}{2}}-\lambda \|\nabla v\|_2$ is strictly increasing in $(0, +\infty)$, and 
$$\lim\limits_{y\to 0}h(y)=-\lambda\|\nabla v\|_2<0,\ \ \ \ \lim\limits_{y\to +\infty}h(y)=+\infty,$$
we see that the equation $h(y)=0$ has a unique solution in $(0, +\infty)$. Therefore, problem (\ref{eq21}) is solvable for $f(x)\in\mathcal{M}$. This completes the proof of Lemma 2.3.


\vskip 0.1in

{\bf Proof of Theorem 2.1:} \ \  To prove Theorem 2.1, we denote by $H^1_0(\Omega)$ the standard Sobolev space with norm $\|u\|=\|\nabla u\|_2$, and consider the following functional defined on $H^1_0(\Omega)$.
$$I_\lambda(u)=\frac{1}{2}\|\nabla u\|^2_2+\frac{b}{2(\alpha+1)}\|\nabla u\|^{2(\alpha+1)}_2-\frac{1}{p+1}\|u^+\|^{p+1}_{p+1}-\lambda\int_\Omega fudx.$$
We claim that $I_\lambda$ is bounded from below on $H^1_0(\Omega)$ and
$$\lim\limits_{\|u\|\rightarrow +\infty} I_\lambda(u)=+\infty.$$

In fact, by H\"{o}lder's and Young's inequality, we get 
\begin{align}
\nonumber
\lambda\int_\Omega fu dx &\leq \frac{\lambda}{\sqrt{\lambda_1}}\|f\|_2\|\nabla u\|_2 \leq \frac{1}{4}\|\nabla u\|^2_2+\frac{\lambda^2}{\lambda_1(\Omega)}\|f\|^2_2
\end{align}
with $\lambda_1(\Omega)$ being the first eigenvalue of Dirichlet Laplacian. By Sobolev's inequality, we have
\begin{align}
\|u^+\|^{p+1}_{p+1}\leq\|u\|^{p+1}_{p+1}\leq S(\Omega)\|\nabla u\|^{p+1}_{2}
\end{align}
for some positive constant $S(\Omega)$ independent of $\lambda$. Therefore,
\begin{equation}\label{21}
\begin{array}{ll}
I_\lambda(u) &\geq\frac{1}{4}\|\nabla u\|^2_2+\frac{b}{2(\alpha+1)}\|\nabla u\|^{2(\alpha+1)}_2-\frac{S(\Omega)}{p+1}\|\nabla u\|^{p+1}_2-\frac{\lambda^2}{\lambda_1(\Omega)}\|f\|^2_2\\
&\geq \frac{b}{2(\alpha+1)}\|\nabla u\|^{2(\alpha+1)}_2-\frac{S(\Omega)}{p+1}\|\nabla u\|^{p+1}_2-\frac{\lambda^2}{\lambda_1(\Omega)}\|f\|^2_2
\end{array}
\end{equation}
which implies that $\lim\limits_{\|u\|\rightarrow +\infty} I_\lambda(u)=+\infty$ due to $1<p<{2\alpha+1}$.

By evaluating the minimum of function $\frac{b}{2(\alpha+1)}t^{2(\alpha+1)}-\frac{S(\Omega)}{p+1}t^{p+1}$ on $(0,+\infty)$, we get
\begin{align}\label{22}
\frac{b}{2(\alpha+1)}\|\nabla u\|^{2(\alpha+1)}_2-\frac{S(\Omega)}{p+1}\|\nabla u\|^{p+1}_2\geq -\frac{\gamma}{2(\alpha+1)(p+1)}[\frac{S^{2(\alpha+1)}(\Omega)}{b^{p+1}}]^{\frac{1}{\gamma}}.
\end{align}
Combining (\ref{21}) and (\ref{22}) together, we have
\begin{equation}\label{23}
I_\lambda(u)\geq -\frac{\gamma}{2(\alpha+1)(p+1)}[\frac{S^{2(\alpha+1)}(\Omega)}{b^{p+1}}]^{\frac{1}{\gamma}}-\frac{\lambda^2}{\lambda_1(\Omega)}\|f\|^2_2.
\end{equation}
This implies that $I_\lambda$ is bounded from below on $H^1_0(\Omega)$.

Setting
\begin{equation}\label{24}
C_\lambda=\inf\limits_{u\in H^1_0(\Omega)} I_\lambda(u),
\end{equation}
we can claim that
\begin{equation}\label{25}
 -\frac{\gamma}{2(\alpha+1)(p+1)}[\frac{S^{2(\alpha+1)}(\Omega)}{b^{p+1}}]^{\frac{1}{\gamma}}-\frac{\lambda^2}{\lambda_1(\Omega)}\|f\|^2_2\leq C_\lambda<0.
\end{equation}
In fact, the first inequality in (\ref{25}) follows from (\ref{23}) and (\ref{24}). To prove the second inequality in (\ref{25}), we denote by $\varphi(x)$ the nontrivial solution of problem (\ref{eq21}). The existence of $\varphi(x)$ follows from Lemma 2.3 since $f(x)\in\mathcal{M}$. Moreover, $\varphi(x)$ verifies
$$\|\nabla\varphi\|_2^2+b\|\nabla\varphi\|_2^{2(\alpha+1)}=\lambda\int_{\Omega}f\varphi dx.$$
Therefore, we have
$$I_\lambda(\varphi)=-\frac{1}{2}\|\nabla\varphi\|_2^2-b(1-\frac{1}{2\alpha+2})\|\nabla\varphi\|_2^{2(\alpha+1)}-\frac{1}{p+1}\|\varphi\|_{p+1}^{p+1}<0.$$
This and the definition of $C_\lambda$ imply
$$C_\lambda\leq I_\lambda(\varphi)<0.$$
By Ekeland's variational principle (see \cite{St}), we know that there exists a sequence $\{u_n\}\subset H^1_0(\Omega)$ such that
\begin{equation}\label{26}
\begin{array}{ll}
\lim\limits_{n\rightarrow +\infty}I_\lambda(u_n)=C_\lambda,\\
\lim\limits_{n\rightarrow +\infty}I'_\lambda(u_n)= 0.
\end{array}
\end{equation}
Since $C_\lambda$ is finite and $\lim\limits_{\|u\|\rightarrow +\infty} I_\lambda(u)=+\infty$, we conclude that $\{u_n\}$ is bounded in $H^1_0(\Omega)$. Therefore, up to a subsequence, we may assume that
\begin{equation}\label{27}\begin{split}
\begin{array}{ll}
u_n\rightharpoonup u\ \ \mbox{weakly in}\ \ H^1_0(\Omega),  \\
u_n\rightarrow u\ \ \mbox{almost everywhere in}\ \ \Omega,   \\
u_n\rightarrow u\ \ \mbox{strongly in}\ \ L^s(\Omega)\ \  \mbox{for any}\ \ s\in (1, 2^*+1),
\end{array}
\end{split}\end{equation}
for some function $u\in H^1_0(\Omega)$.

Consequently, we have
\begin{equation}\label{281}
\begin{array}{ll}
\int_\Omega((u^+_n)^p+\lambda f)(u_n-u)dx\rightarrow 0,\\
\int_\Omega\nabla u_n\cdot\nabla u dx\rightarrow\int_\Omega |\nabla u|^2 dx.
\end{array}
\end{equation}
Since
$$\langle I'_\lambda(u_n),u_n-u\rangle=(1+b\|\nabla u_n\|^{2\alpha}_2)\int_\Omega\nabla u_n\cdot\nabla(u_n-u)dx-\int_\Omega((u^+_n)^p+\lambda f)(u_n-u)dx,$$
it follows from (\ref{281}) and the fact $\langle I'_\lambda(u_n),u_n-u\rangle\rightarrow 0$ that
$$(1+b\|\nabla u_n\|^{2\alpha}_2)\int_\Omega\nabla u_n\cdot\nabla(u_n-u)dx\rightarrow0.$$
This implies $\|u_n\|\rightarrow \|u\|$. Therefore, $u_n\rightarrow u$ strongly in $H^1_0(\Omega)$.

For any $\phi\in H^1_0(\Omega)$, we have
$$\langle I'_\lambda(u_n),\phi\rangle=(1+b\|\nabla u_n\|^{2\alpha}_2)\int_\Omega\nabla u_n\cdot\nabla\phi dx-\int_\Omega((u^+_n)^p+\lambda f)\phi dx.$$
By sending $n$ to $+\infty$ in the above equation, we get
$$(1+b\|\nabla u\|^{2\alpha}_2)\int_\Omega\nabla u\cdot\nabla\phi dx=\int_\Omega(u^+)^p\phi+\lambda\int_\Omega f\phi dx.$$
Therefore, $u$ is a weak solution of the following problem
\begin{equation}\label{eq28}\begin{split}
\left\{\begin{array}{ll}
-(1+b\|\nabla u\|^{2\alpha}_2)\Delta u=(u^+)^p+\lambda f(x) & x\in\Omega, \\
u=0 & x\in\partial\Omega.
\end{array}
\right.
\end{split}\end{equation}
Furthermore, we can prove $u(x)$ is positive in $\Omega$ by strong comparison principle of Laplace operator. In fact, by the assumption $f(x)\in\mathcal{M}$,  we know from Lemma 2.3 that there exists a function $\varphi(x)$ which satisfies. 
\begin{equation}\label{28}
\left \{ \begin{array}{ll}
-(1+b\|\nabla\varphi\|_2^{2\alpha)})\Delta\varphi=\lambda f(x) & x\in\Omega,\\
\varphi\geq 0 & x\in\Omega,\\
\varphi=0 & x\in\partial\Omega.
\end{array}
\right.
\end{equation}
By (\ref{eq28}) and (\ref{28}), we can easily see 
\begin{equation}\label{29}
\left\{\begin{array}{ll}
-\frac{1+b\|\nabla\varphi\|_2^{2\alpha}}{1+b\|\nabla u\|_2^{2\alpha}}\Delta\varphi\leq -\Delta u & x\in\Omega,\\
\varphi=u=0 & x\in\partial\Omega.
\end{array}
\right.
\end{equation}
Therefore, by comparison principle for weak solutions, we have
$$u(x)\geq\frac{1+b\|\nabla\varphi\|_2^{2\alpha}}{1+b\|\nabla u\|_2^{2\alpha}}\varphi(x)\geq 0\ \ \mbox{for any}\ \ x\in\Omega.$$
This and (\ref{eq28}) imply that $u$ is a nonnegative weak solution of problem (\ref{eq3}). Moreover, by regularity theory of elliptic equations, we know further that $u$ is a nonnegative classical solution of (\ref{eq3}). Finally, by strong comparison principle of Laplace operator, we have
$$u(x)>\frac{1+b\|\nabla\varphi\|_2^{2\alpha}}{1+b\|\nabla u\|_2^{2\alpha}}\varphi(x)\geq 0\ \ \mbox{for any}\ \ x\in\Omega.$$
Therefore, $u$ is a positive solution of problem (\ref{eq3}), and the proof of Theorem 2.1 is completed.

To prove Theorem 2.2, we need the following result which was proven in \cite{DLS}.
\vskip 0.1in

{\em {\bf Lemma 2.4}(\cite{DLS}) If $1<p<2\alpha+1$ and $b>b_0$, then the following problem has no solution.
\begin{equation}\label{eq29}\begin{split}
\left\{\begin{array}{ll}
-(1+b\|\nabla u\|^{2\alpha}_2)\Delta u=u^p & x\in\Omega, \\
u>0 & x\in\Omega,\\
u=0 & x\in\partial\Omega.
\end{array}
\right.
\end{split}
\end{equation} }

The following lemma is crucial for proving Theorem 2.2.
\vskip 0.1in
{\em {\bf Lemma 2.5} If $1<p<2\alpha+1$, $b>b_0$ and $u_\lambda(x)$ is a positive solution of problem (\ref{eq3}) corresponding to parameter $\lambda$, then we have 
$$\|u_\lambda\|_{\infty}\rightarrow 0,  \  \mbox{ as}\    \lambda\rightarrow 0.$$}
\vskip 0.1in

{\bf Proof:}\ \ We adopt a contradiction argument. Suppose that the conclusion of Lemma 2.5 is not true, then there would exist a sequence $\{\lambda_n\}^\infty_{n=1}\subset (0, 1)$, and $\{u_{\lambda_n}\}^\infty_{n=1}$ such that

\begin{equation}\label{eq210}\begin{split}
\left\{\begin{array}{ll}
-(1+b\|\nabla u_{\lambda_n}\|^{2\alpha}_2)\Delta u_{\lambda_n}=u^p_{\lambda_n}+\lambda_n f(x) & x\in\Omega, \\
u_{\lambda_n}\geq 0 & x\in\Omega,\\
u_{\lambda_n}=0 & x\in\partial\Omega,
\end{array}
\right.
\end{split}\end{equation}
and $M_n=\|u_{\lambda_n}\|_{\infty}\rightarrow C>0, \lambda_n\rightarrow 0,  as\   n\rightarrow+\infty.$

Since $1<p<2\alpha+1$, we get easily from (\ref{eq210}) that
$$\|\nabla u_{\lambda_n}\|_\infty\leq C,$$
for some positive constant $C$ independent of $n$. Furthermore, by a bootstrap argument and Schauder's estimates of elliptic equations, we have
$$\|u_{\lambda_n}\|_{C^{2,\tau}(\Omega)}\leq C_1,$$
for some constant $C_1$ independent of $n$ and $\tau\in(0,1)$. Therefore, up to a subsequence, $u_{\lambda_n}$ converges in $C^2(\Omega)$ to a nonnegative function $u$ which satisfies

\begin{equation}\label{eq211}\begin{split}
\left\{\begin{array}{ll}
-(1+b\|\nabla u\|^{2\alpha}_2)\Delta u=u^p  & x\in\Omega, \\
u\geq 0 & x\in\Omega,\\
u=0 & x\in\partial\Omega.
\end{array}
\right.
\end{split}\end{equation}

Since $\|u\|_{\infty}=\lim\limits_{n\rightarrow +\infty}\|u_{\lambda_n}\|_{\infty}=C>0$, we can deduce from the strong maximum principle that $u(x)>0$ for any $x\in\Omega$. Therefore, $u(x)$ is a solution of problem (\ref{eq29}). This contradicts Lemma 2.4. 

\vskip 0.1in

{\em {\bf Lemma 2.6} If $1<p<2\alpha+1$, $\alpha\geq\frac{1}{2}$ and $b>b_0$, then problem (\ref{eq3}) has at most one positive solution for parameter $\lambda$ small enough.}
\vskip 0.1in
{\bf Proof:}\ \  Let $u_\lambda(x)$ and $v_\lambda(x)$ be two arbitrary positive solutions of problem (\ref{eq3}). That is, $u_\lambda(x)$ and $v_\lambda(x)$ satisfy

\begin{equation}\label{eq212}
\begin{split}
\left\{
\begin{array}{ll}
-(1+b\|\nabla u_{\lambda}\|^{2\alpha}_2)\Delta u_{\lambda}=u^p_{\lambda}+\lambda f(x) & x\in\Omega, \\
-(1+b\|\nabla v_{\lambda}\|^{2\alpha}_2)\Delta v_{\lambda}=v^p_{\lambda}+\lambda f(x) & x\in\Omega, \\
u_{\lambda}=v_\lambda=0 & x\in\partial\Omega.
\end{array}
\right.
\end{split}
\end{equation}
What we should do is that $u_\lambda(x)\equiv v_\lambda(x)$ in $\Omega$ for small enough parameter $\lambda$. To this end, we set $w_\lambda(x)=u_\lambda(x)-v_\lambda(x)$, and $A=b(\|\nabla u_\lambda\|_2^{2\alpha}-\|\nabla v_\lambda\|_2^{2\alpha})$. By (\ref{eq212}) and mean value theorem, we know that there exists a function $0\leq\theta(x)\leq 1$  such that $w_\lambda(x)$ verifies
$$-(1+b\|\nabla u_{\lambda}\|^{2\alpha}_2)\Delta w_{\lambda}=p(\theta u_{\lambda}+(1-\theta)v_\lambda)^{p-1}w_\lambda+A\Delta v_\lambda\ \ \ \ \ \ x\in\Omega.$$
Multiplying the above equality by $w_\lambda$ and integrating on $\Omega$, we get
\begin{equation}\label{eq213}
(1+b\|\nabla u_{\lambda}\|^{2\alpha}_2)\|\nabla w_{\lambda}\|_2^2=p\int_\Omega(\theta u_{\lambda}+(1-\theta)v_\lambda)^{p-1}w^2_\lambda dx-A\int_\Omega\nabla v_\lambda\centerdot\nabla w_\lambda dx.
\end{equation}
By mean value theorem and triangle inequality, we have
\begin{equation}\label{eq214}
\begin{array}{ll}
|-A|&=2\alpha b|(\theta_0\|\nabla u_\lambda\|_2+(1-\theta_0)\|\nabla v_\lambda\|_2)^{2\alpha-1}(\|\nabla u_\lambda\|_2-\|\nabla v_\lambda\|_2)|\\
&\leq 2\alpha b(\|\nabla u_\lambda\|_2+\|\nabla v_\lambda\|_2)^{2\alpha-1}\|\nabla w_\lambda\|_2\\
&\equiv C_1(\lambda)\|\nabla w_\lambda\|_2.
\end{array}
\end{equation}
Where $C_1(\lambda)=2\alpha b(\|\nabla u_\lambda\|_2+\|\nabla v_\lambda\|_2)^{2\alpha-1}$.

Since $(\theta u_{\lambda}+(1-\theta)v_\lambda)^{p-1}\leq (\|u_{\lambda}\|_{L^\infty(\Omega)}+\|v_{\lambda}\|_{L^\infty(\Omega)})^{p-1}$,  by Poincare inequality we have
\begin{equation}\label{eq215}
\begin{array}{ll}
|p\int_\Omega(\theta u_{\lambda}+(1-\theta)v_\lambda)^{p-1}w^2_\lambda dx|&\leq p(\|u_{\lambda}\|_{L^\infty(\Omega)}+\|v_{\lambda}\|_{L^\infty(\Omega)})^{p-1}\int_\Omega w_\lambda^2 dx\\
&\leq\frac{p}{\lambda_1(\Omega)}|(\|u_{\lambda}\|_{L^\infty(\Omega)}+\|v_{\lambda}\|_{L^\infty(\Omega)})^{p-1}\|\nabla w_\lambda\|_2^2\\
&\equiv C_2(\lambda)\|\nabla w_\lambda\|_2^2.
\end{array}
\end{equation}
Where $C_2(\lambda)=\frac{p}{\lambda_1(\Omega)}|(\|u_{\lambda}\|_{L^\infty(\Omega)}+\|v_{\lambda}\|_{L^\infty(\Omega)})^{p-1}$ and $\lambda_1(\Omega)$ is the first eigenvalue of the Dirichlet Laplacian.

From (\ref{eq213}), (\ref{eq214}) and (\ref{eq215}), we get
\begin{equation}\label{eq216}
\|\nabla w_\lambda\|_2^2\leq (C_2(\lambda)+\|\nabla v_\lambda\|_2C_1(\lambda))\|\nabla w_\lambda\|_2^2.
\end{equation}
Since $p>1$ and $2\alpha-1\geq 0$, by Lemma 2.5 we know that
\begin{equation}\label{eq217}
\lim\limits_{\lambda\rightarrow 0}(C_2(\lambda)+\|\nabla v_\lambda\|_2C_1(\lambda))=0.
\end{equation}

Combining (\ref{eq216}) and (\ref{eq217}) together imply that there exists a positive number $\lambda_0$ such that $\|\nabla w_\lambda\|_2=0$ for any $\lambda\in (0, \lambda_0)$. Therefore, $w_\lambda\equiv  0$ in $\Omega$ for any $\lambda\in (0, \lambda_0)$ because $w_\lambda=0$ on $\partial\Omega$. This completes the proof of Lemma 2.6.
\vskip 0.1in

{\bf Proof of Theorem 2.2:}\ \  In the sequel, we always assume that $1<p<2\alpha+1$ and $b>b_0$. If $f(x)\in\mathcal{M}$, then Theorem 2.1 guarantees the  existence of positive solution for problem (\ref{eq3}). If in addition $\alpha\geq\frac{1}{2}$, then Lemma 2.6 implies that the uniqueness claim in Theorem 2.2 is true. Therefore, to complete the proof of Theorem 2.2, we just need to prove that the necessary condition for positive solvability of problem (\ref{eq3}) for $\lambda>0$ is $f(x)\in\mathcal{M}$. To make this end, we assume that  problem (\ref{eq3}) has positive solution for any $\lambda>0$. Let $u_\lambda$ be positive solution of problem (\ref{eq3}) with respect to parameter $\lambda$. By Lemma 2.5, we have

$$\|u_\lambda\|_{\infty}\rightarrow 0\ \  \mbox{as}\ \  \lambda\rightarrow 0.$$

 Let $u_\lambda=\lambda v_\lambda$, then $v_\lambda$ satisfies
\begin{equation}\label{eq218}\begin{split}
\left\{\begin{array}{ll}
-(1+b\lambda^{2\alpha}\|\nabla v_\lambda\|^{2\alpha}_2)\Delta v_\lambda=\lambda^{p-1}v^p_\lambda+f(x)  & \mbox{in}\ \ \Omega, \\
v_\lambda\geq0& \mbox{in}\ \ \Omega, \\
v_\lambda=0 & \mbox{on}\ \ \partial\Omega.
\end{array}
\right.
\end{split}\end{equation}
Multiplying the differential equation in problem (\ref{eq218}) by $v_\lambda$ and integrating the result equation over $\Omega$, we get
\begin{equation}\label{eq219}
(1+b\lambda^{2\alpha}\|\nabla v_\lambda\|^{2\alpha}_2)\|\nabla v_\lambda\|^2_2=\int_\Omega \lambda^{p-1}v^{p+1}_\lambda dx+\int_\Omega fv_\lambda dx,
\end{equation}
that is,
\begin{equation}\label{eq220}
(1+b\lambda^{2\alpha}\|\nabla v_\lambda\|^{2\alpha}_2)\|\nabla v_\lambda\|^2_2=\int_\Omega u_\lambda^{p-1}v^2_\lambda dx+\int_\Omega fv_\lambda dx.
\end{equation}

Denote by $\lambda_1(\Omega)$ the first eigenvalue of Dirichlet Laplacian. By H\"{o}lder's, Poincare's and Young's inequality, we have
\begin{equation}\label{eq221}
|\int_\Omega fv_\lambda dx|\leq \|f\|_2\|v_\lambda\|_2\leq \frac{1}{4}\|\nabla v_\lambda\|^2_2+\frac{\|f\|^2_2}{\lambda_1(\Omega)}.
\end{equation}
 Since $\lim\limits_{\lambda\to 0}\|u_\lambda\|_{\infty}=0$, there is a positive constant $\lambda_0$ such that
\begin{equation}
\|u_\lambda\|_{\infty}\leq (\frac{\lambda_1(\Omega)}{4})^{\frac{1}{p-1}},\  \ \ \ \ \mbox{for}\ \ \ \  \lambda\in(0,\lambda_0).
\end{equation}
From this and Poincare's inequality, we have
\begin{equation}\label{eq222}
\int_\Omega u^{p-1}_\lambda v^2_\lambda dx\leq\|u_\lambda\|^{p-1}_{\infty}\|v_\lambda\|_2^2\leq\frac{1}{4}\|\nabla v_\lambda\|^2_2\ \ \ \ \ \mbox{for}\ \ \ \ \ \lambda\in(0, \lambda_0).
\end{equation}
Combining (\ref{eq220}), (\ref{eq221}) and (\ref{eq222}) together, we get
$$\|\nabla v_\lambda\|^2_2\leq \frac{2\|f\|^2_2}{\lambda_1(\Omega)}\  \ \ \ \  \mbox{for}\ \ \\ \lambda\in(0,\lambda_0).$$
The above inequality and a bootstrap argument show that there exists a positive constant $C$ independent of $\lambda$ such that
$$\|v_\lambda\|_{\infty}\leq C\ \ \ \    \mbox{for}\  \lambda\in(0,\lambda_0).$$
Furthermore, by standard elliptic regularity theory, we can find a positive constant $C$ independent of $\lambda$ such that
$$\|v_\lambda\|_{C^{2,\tau}(\Omega)}\leq C\ \ \ \ \ \mbox{for some}\ \tau\in(0, 1)\ \mbox{and any}\ \lambda\in(0,\lambda_0).$$
Therefore, up to a subsequence, we may assume that
$$v_\lambda\rightarrow v\geq 0\ \ \ \ \mbox{in $C^2(\Omega)$} \ \mbox{as}\ \lambda\rightarrow 0.$$

Sending $\lambda$ to $0$ in problem (\ref{eq218}), we see that $v$ verifies

\begin{equation}\label{eq223}\begin{split}
\left\{\begin{array}{ll}
-\Delta v=f(x) & \mbox{in}\ \ \Omega, \\
v\geq 0 & \mbox{in}\ \ \Omega, \\
v=0 & \mbox{on}\ \ \partial\Omega.
\end{array}
\right.
\end{split}\end{equation}
Therefore,  $f(x)\in\mathcal{M}$. This completes the proof of Theorem 2.2.

\section*{3. The case $2\alpha+1<p<2^*$}

\setcounter{section}{3}

\setcounter{equation}{0}

\noindent

This section devotes to deal with the case $2\alpha+1<p<2^*$. The main purpose is to prove the following result.
\vskip 0.1in
{\em {\bf Theorem 3.1}\ \ If $2\alpha+1<p<2^*$ and $f(x)\in\mathcal{M}$, then there are two positive constants $\lambda_{pf}, \Lambda_{pf}<+\infty$ such that problem (\ref{eq3}) has at least two positive solutions for $\lambda\in (0, \lambda_{pf})$, and has no positive solution for $\lambda>\Lambda_{pf}$.} 
\vskip 0.1in
{\em {\bf Remark 3.2}\ \ Instead of multiplicity results, if we focus only on the existence result, then the condition $f(x)\in\mathcal{M}$ may be made a small relaxation (see Lemma 3.3 of this section).}

To prove Theorem 3.1, we denote by $H_0^1(\Omega)$ the standard Sobolev space, and consider functional
\begin{equation}\label{eq37}
J_\lambda(u)=\frac{1}{2}\|\nabla u\|_2^2+\frac{b}{2(\alpha+1)}\|\nabla u\|_2^{2(\alpha+1)}-\frac{1}{p+1}\|u^+\|^{p+1}_{p+1} dx-\lambda\int_\Omega fu dx
\end{equation}
defined on $H_0^1(\Omega)$. It is obvious that any critical point $u\in H_0^1(\Omega)$ of $J_\lambda(u)$ is a weak solution of problem
\begin{equation}\label{eq38}
\left\{
\begin{array}{ll}
-(1+b\|\nabla u\|_2^{2\alpha})\Delta u=(u^+)^p+\lambda f & x\in\Omega,\\
u=0 & x\in\partial\Omega.
\end{array}
\right.
\end{equation}
Let $\mathcal{N}(\partial\Omega)\subset\overline{\Omega}$ denote inner neighborhood of $\partial\Omega$. Setting
$$\mathcal{F}^+=\{\ f\in C^1(\overline{\Omega}\backslash\{0\}\ \mbox{with property that}\ f(x)\geq 0\ \mbox{in some}\ \mathcal{N}(\partial\Omega)\ \}.$$
Obviously, $\mathcal{F}^+\backslash\mathcal{M}\neq\O$. In fact, any nontrivial function $\phi(x)$ with property $\phi(x)\leq 0$ in $\Omega$ and $supp\{\phi(x)\}\subset\Omega$ belongs to $\mathcal{F}^+$, but not belongs to $\mathcal{M}$. Instead of condition $f(x)\in\mathcal{M}$, we will find a positive solution of problem (\ref{eq3}) in the following lemma under the condition $f(x)\in\mathcal{M}\cup\mathcal{F}^+$.
\vskip 0.1in
{\em {\bf Lemma 3.3}\ \ If $2\alpha+1<p<2^*$ and $f(x)\in\mathcal{M}\cup\mathcal{F}^+$, then there exists a positive number $\lambda_f$ such that problem (\ref{eq3}) has a positive solution $v_\lambda$ for any $\lambda\in (0, \lambda_f)$ with property that $J_\lambda(v_\lambda)>0$ and $v_\lambda$ converges, as $\lambda\rightarrow 0$, to a solution $v$ of the following problem
\begin{equation}\label{eq39}
\left\{
\begin{array}{ll}
-(1+b\|\nabla v\|_2^{2\alpha})\Delta v=v^p & x\in\Omega,\\
v>0 & x\in\Omega,\\
v=0 & x\in\partial\Omega.
\end{array}
\right.
\end{equation}
}
\vskip 0.1in

{\bf Proof:}\ \ We prove this lemma by the following steps.
\vskip 0.1in

{\it Step1:}\ \ There are positive numbers $\beta_f$, $\rho_0$, $E_0$ and elements $e_0, e_1\in H_0^1(\Omega)$ independent of $\lambda$ such that
$$\|\nabla e_0\|_2<\rho_0<\|\nabla e_1\|_2\ \ \mbox{and}\ \ J_\lambda(u)|_{\partial B_{\rho_0}}\geq E_0>\max\{J_\lambda(e_0), J_\lambda(e_1)\}$$
for any $\lambda\in (0, \beta_f)$. Where $B_{\rho_0}=\{u\in H_0^1(\Omega):\ \|\nabla u\|_2<\rho_0\}$.

In fact, if we denote by $\lambda_1(\Omega)$ the first eigenvalue of the eigenvalue problem
\begin{equation}\label{eq310}\begin{split}
\left\{\begin{array}{ll}
-\Delta \phi=\lambda\phi, & x\in\Omega, \\
\phi=0, & x\in\partial\Omega,
\end{array}
\right.
\end{split}\end{equation}
then we have
\begin{align}\label{eq311}
|\lambda\int_\Omega fu dx| &\leq \frac{\lambda}{\sqrt{\lambda_1(\Omega)}}\|f\|_2\|\nabla u\|_2 \\
&\leq \frac{1}{4}\|\nabla u\|^2_2+\frac{\lambda^2}{\lambda_1(\Omega)}\|f\|^2_2.
\end{align}
Therefore,
\begin{equation}\label{eq312}
J_\lambda(u)\geq\frac{1}{4}\|\nabla u\|^2_2-\frac{1}{p+1}\|u\|^{p+1}_{p+1}-\frac{\lambda^2}{\lambda_1(\Omega)}\|f\|^2_2.
\end{equation}
By the assumption $1<p<2^*$ and Sobolev's inequality, we have
\begin{equation}\label{eq313}
\|u^+\|^{p+1}_{p+1}\leq\|u\|^{p+1}_{p+1}\leq S(\Omega)\|\nabla u\|^{p+1}_{2}
\end{equation}
with $S(\Omega)$ independent of $\lambda$.

Combining (\ref{eq312}) and (\ref{eq313}) together, we get
$$J_\lambda(u)\geq\frac{1}{4}\|\nabla u\|^2_2-\frac{S(\Omega)}{p+1}\|\nabla u\|^{p+1}_2-\frac{\lambda^2}{\lambda_1(\Omega)}\|f\|^2_2.$$
Hence
$$J_\lambda(u)|_{\partial B_\rho}\geq\frac{1}{4}\rho^2-\frac{S(\Omega)}{p+1}\rho^{p+1}-\frac{\lambda^2}{\lambda_1(\Omega)}\|f\|^2_2.$$
Noting $p+1>2$, we can choose positive number $\rho_0$ independent of $\lambda$ so small that
$$\frac{1}{4}\rho^2_0-\frac{S(\Omega)}{p+1}\rho_0^{p+1}=E_1>0.$$
Taking
$$\beta_f=\frac{\sqrt{3\alpha_1\lambda_1(\Omega)}}{2\|f\|_2}\ \mbox{and}\ \ E_0=\frac{E_1}{4},$$
we get
\begin{equation}\label{eq314}
J_\lambda(u)|_{\partial B_{\rho_0}}\geq E_0>0,\   \mbox{for}\  \lambda\in(0,\beta_f).
\end{equation}
Since $J_\lambda(0)=0$, we may take $e_0=0$. To choose a suitable $e_1$, we denote by $\phi_1(x)$ the first eigenfunction corresponding to $\lambda_1(\Omega)$. By the definition of $J_\lambda(u)$, we have
$$J_\lambda(t\phi_1)\leq\frac{\|\nabla \phi_1\|^2_2}{2}t^2+\frac{b\|\nabla\phi_1\|^{2(\alpha+1)}_2}{2(\alpha+1)}t^{2(\alpha+1)}-\frac{\|\phi_1\|^{p+1}_{p+1}}{p+1}t^{p+1}+\beta_f\|f\|_2\|\phi_1\|_2t$$
for any $\lambda\in (0, \beta_f)$. Noting $p+1>2(\alpha+1)>2$, we have
$$\lim\limits_{t\rightarrow +\infty}(\frac{\|\nabla \phi_1\|^2_2}{2}t^2+\frac{b\|\nabla\phi_1\|^{2(\alpha+1)}_2}{2(\alpha+1)}t^{2(\alpha+1)}-\frac{\|\phi_1\|^{p+1}_{p+1}}{p+1}t^{p+1}+\beta_f\|f\|_2\|\phi_1\|_2t)=-\infty.$$
Therefore, we can choose a large constant $t_0$ independent of $\lambda$ such that $t_0\|\nabla\phi_1\|_2>\rho_0$ and
$$\frac{\|\nabla \phi_1\|^2_2}{2}t_0^2+\frac{b\|\nabla\phi_1\|^{2(\alpha+1)}_2}{2(\alpha+1)}t_0^{2(\alpha+1)}-\frac{\|\phi_1\|^{p+1}_{p+1}}{p+1}t_0^{p+1}+\beta_f\|f\|_2\|\phi_1\|_2t_0<0.$$ 
Taking $e_1=t_0\phi_1(x)$,  we have $J_\lambda(e_1)=J_\lambda(t_0\phi_1)<0$ for any $\lambda\in (0, \beta_f)$.  In summary, for the above choices of $\beta_f$, $\rho_0$, $E_0$, $e_0$ and $e_1$, we have
$$\|\nabla e_0\|_2<\rho_0<\|\nabla e_1\|_2\ \ \mbox{and}\ \ J_\lambda(u)|_{\partial B_{\rho_0}}\geq E_0>\max\{J_\lambda(e_0), J_\lambda(e_1)\}$$
for any $\lambda\in (0, \beta_f)$. This concludes {\it Step1}.
\vskip 0.1in

{\it Step2:}\ \ For any $\lambda\in (0, \beta_f)$, problem (\ref{eq38}) has a solution $v_\lambda(x)$ with property $J_\lambda(v_\lambda)\geq E_0>0$.

To conclude {\it Step 2}, for any $\lambda\in (0, \beta_f)$, we set
$$\Gamma=\{\gamma\in C([0,1],H^1_0(\Omega)):\ \gamma(0)=e_0=0,\ \gamma(1)=e_1=t_0\phi_1\},$$
and
$$C_\lambda=\inf\limits_{\gamma\in\Gamma}\max\limits_{s\in[0,1]}J_\lambda(\gamma(s)).$$
where $t_0$ is a fixed constant given in {\it Step1}.

By {\it Step1} and mountain pass theorem without $PS$ condition, we know that there is a sequence $\{v^n_\lambda\}\subset H^1_0(\Omega)$ such that
\begin{equation}\label{eq315}
\begin{array}{ll}
J_\lambda(v^n_\lambda)\rightarrow C_\lambda\geq E_0>0 &\mbox{as}\  n\rightarrow+\infty,\\
J'_\lambda(v^n_\lambda)\rightarrow0 &\mbox{as}\  n\rightarrow+\infty.
\end{array}
\end{equation}

Because of $2\alpha+1<p<2^*$, it is easy to verify that $J_\lambda(u)$ satisfies $PS$ condition. Therefore, up to a subsequence, $\{v^n_\lambda\}$ converges strongly in $H^1_0(\Omega)$ to a function $v_\lambda$ which satisfies

\begin{equation}\label{eq316}\begin{split}
\left\{\begin{array}{ll}
-(1+b\|\nabla v_\lambda\|^{2\alpha}_2)\Delta v_\lambda=(v^+_\lambda)^p+\lambda f(x)& \mbox{in}\ \ \Omega, \\
v_\lambda=0 & \mbox{on}\ \ \partial\Omega.
\end{array}
\right.
\end{split}\end{equation}
and
\begin{equation}\label{eq317}
\left\{\begin{array}{ll}
J_\lambda(v_\lambda)=C_\lambda\geq E_0>0,\\
J^\prime_\lambda(v_\lambda)=0.
\end{array}
\right.
\end{equation}
This makes {\it Step2}.
\vskip 0.1in

{\it Step3:}\ \  There exists a positive number $\lambda_f\leq\beta_f$ such that, for any $\lambda\in (0, \lambda_f)$, the solution $v_\lambda(x)$ obtained in {\it Step 2} for problem (\ref{eq38}) is positive. Therefore, $v_\lambda(x)$ is a positive solution to problem (\ref{eq3}) and $J_\lambda(v_\lambda)\geq E_0>0$ for any $\lambda\in (0, \lambda_f)$.

We divide the proof of {\it Step3} into two cases. One is $f(x)\in\mathcal{M}$,  and the other is $f(x)\in\mathcal{F}^+$.

If $f(x)\in\mathcal{M}$, then by Lemma 2.3, we know that problem (\ref{eq21}) has a nonnegative solution $u_{0,\lambda}(x)$ for any $\lambda\in (0, \beta_f)$. Since $v_\lambda(x)$ is a solution of (\ref{eq38}) for $\lambda\in (0, \beta_f)$, we have
\begin{equation}\label{eq318}
\left\{\begin{array}{ll}
-(1+b\|\nabla v_\lambda\|_2^{2\alpha})\Delta v_\lambda\geq\lambda f=-(1+b\|\nabla u_{0,\lambda}\|_2^{2\alpha})\Delta u_{0,\lambda} & x\in\Omega,\\
v_\lambda=u_{0,\lambda}=0 &x\in\Omega,
\end{array}
\right.
\end{equation}
for any $\lambda\in (0, \beta_f)$. Therefore, by strong comparison principle of Laplace operator, we have 
$$v_\lambda(x)>\frac{1+b\|\nabla u_{0,\lambda}\|_2^{2\alpha}}{1+b\|\nabla v_\lambda\|_2^{2\alpha}} u_{0,\lambda}(x)\geq 0\ \ \ x\in\Omega$$
for any $\lambda\in (0, \beta_f)$.

If $f(x)\in\mathcal{F}^+$, we first claim that there exists a positive constant $C$ independent of $\lambda$ such that
\begin{equation}\label{eq318}
0<E_0\leq C_\lambda\leq C, \ \ \mbox{for}\ \lambda\in (0, \beta_f)
\end{equation}
where $C_\lambda$ is the critical value defined in {\it Step2}, and $\alpha_0$ is the constant given in {\it Step1}. 

In fact, for any $\gamma(s)\in\Gamma$, $g(s)=\|\nabla\gamma(s)\|_2$ is continuous in $[0, 1]$. Since $0=g(0)<\rho_0<\|\nabla e_1\|=g(1)$, by intermediate value theorem, we have $g(s_0)=\rho_0$ for some $s_0\in (0, 1)$. Hence, for any $\gamma(s)\in\Gamma$, we can conclude from {\it Step1} that
$$\max\limits_{s\in[0,1]}J_\lambda(\gamma(s))\geq J_\lambda(\gamma(s_0))\geq E_0.$$
Therefore, for any $\lambda\in (0, \beta_f)$, there holds
$$C_\lambda=\inf\limits_{\gamma\in\Gamma}\max\limits_{s\in[0,1]}J_\lambda(\gamma(s))\geq E_0>0.$$
To derive a upper bound of $C_\lambda$, we take $\gamma_0(s)=se_1=st_0\phi_1$. Obviously, $\gamma_0(s)\in\Gamma$. By the definition of $C_\lambda$, we have
$$C_\lambda\leq\max\limits_{s\in[0, 1]}J_\lambda(se_1)=\max\limits_{t\in[0, t_0]}J_\lambda(t\phi_1).$$
For $t\in[0, t_0]$ and $\lambda\in (0, \beta_f)$, we can get from the definition of $J_\lambda(u)$ that
$$J_\lambda(t\phi_1)\leq\frac{\|\nabla \phi_1\|^2_2}{2}t_0^2+\frac{b\|\nabla\phi_1\|^{2(\alpha+1)}_2}{2(\alpha+1)}t_0^{2(\alpha+1)}+\frac{\|\phi_1\|^{p+1}_{p+1}}{p+1}t_0^{p+1}+\beta_f\|f\|_2\|\phi_1\|_2t_0.$$
Setting
$$C=\frac{\|\nabla \phi_1\|^2_2}{2}t_0^2+\frac{b\|\nabla\phi_1\|^{2(\alpha+1)}_2}{2(\alpha+1)}t_0^{2(\alpha+1)}+\frac{\|\phi_1\|^{p+1}_{p+1}}{p+1}t_0^{p+1}+\beta_f\|f\|_2\|\phi_1\|_2t_0,$$
we see that $C$ is independent of $\lambda$, and
$$C_\lambda\leq C\ \ \mbox{for}\ \lambda\in (0, \beta_f).$$
Therefore, claim (\ref{eq318}) is valid
 
Taking $1<2\alpha+1<p<2^*$ into account, we can conclude from (\ref{eq317}) and (\ref{eq318}) that there exists a positive constant $C$ independent of $\lambda$ such that
$$\|\nabla v_\lambda\|_2\leq C\ \ \mbox{for}\ \ \lambda\in (0, \beta_f).$$
By bootstrap argument and standard regularity theory of elliptic equations, we can conclude from the above estimate that
\begin{equation}\label{eq319}
\|v_\lambda\|_{C^{2,\tau}(\Omega)}\leq C
\end{equation}
for $\lambda\in (0, \beta_f)$, some positive constant $C$ independent of $\lambda$, and $\tau\in (0, 1)$.

Next, we show that $v_\lambda$ is positive in $\Omega$. Since $f(x)\in\mathcal{F}^+$, there exists a neighborhood $\mathcal{N}(\partial\Omega)$ of $\partial\Omega$ such that $f(x)\geq 0$ for $x\in\mathcal{N}(\partial\Omega)$. Set $\Omega_0=\overline{\Omega\backslash\mathcal{N}(\partial\Omega)}$. At first, we can claim that there exists a positive constant $\lambda_f\leq\beta_f$ such that
$$v_\lambda(x)>0,\    \mbox{for}\    x\in\Omega_0, \lambda\in(0,\lambda_f).$$
Otherwise, there would exist a sequence $\lambda_n\rightarrow 0$ as $n\rightarrow +\infty$, and a sequence $x_n\in\Omega_0$ such that
\begin{equation}\label{eq320}
v_{\lambda_n}(x_n)\leq 0,\    \mbox{for}\  n=1,2,\cdots
\end{equation}
By (\ref{eq319}), up to a subsequence, we may assume that $\{v_{\lambda_n}\}$ converges in $C^2(\Omega)$ to function $v$ which satisfies
\begin{equation}\label{eq321}\begin{split}
\left\{\begin{array}{ll}
-(1+b\|\nabla v\|^{2\alpha}_2)\Delta v=(v^+)^p &  x\in\Omega, \\
v=0 & x\in\partial\Omega.
\end{array}
\right.
\end{split}\end{equation}
Noticing that
$$\frac{1}{2}\|\nabla v\|^2_2+\frac{b}{2(\alpha+1)}\|\nabla v\|^{2(\alpha+1)}_2-\frac{1}{p+1}\int_\Omega(v^+)^{p+1}dx=\lim\limits_{n\rightarrow+\infty}C_{\lambda_n}\geq E_0>0,$$
we have $v\not\equiv 0$. Therefore, by strong maximum principle, we have
$$v(x)>0,\  \mbox{for}\    x\in\Omega.$$
In particular,
\begin{equation}\label{eq322}
v(x)>0,\     \mbox{for}\    x\in\Omega_0.
\end{equation}
Because $\Omega_0$ is closed and bounded, we may assume that $\lim\limits_{n\rightarrow\infty}x_n=x_0\in\Omega_0$. Consequently, by (\ref{eq320}), we have  
$$v(x_0)=\lim\limits_{n\rightarrow\infty}v_{\lambda_n}(x_n)\leq 0.$$
This contradicts (\ref{eq322}).

On the second, we can easily see that $v_\lambda(x)>0$ in $\mathcal{N}(\partial\Omega)$ for $\lambda\in(0,\lambda_f)$. In fact, for any $\lambda\in (0, \lambda_f)$, $v_\lambda$ satisfies
\begin{equation}\label{eq323}\begin{split}
\left\{\begin{array}{ll}
-(1+b\|\nabla v_\lambda\|^{2\alpha}_2)\Delta v_\lambda=(v^+_\lambda)^p+\lambda f(x) & \mbox{in}\ \ \Omega, \\
u_\lambda=0 & \mbox{on}\ \ \partial\Omega.
\end{array}
\right.
\end{split}\end{equation}
Therefore, for any $\lambda\in (0, \lambda_f)$, we have 
$$-(1+b\|\nabla v_\lambda\|^{2\alpha}_2)\Delta v_\lambda\geq 0,\ \ \ \mbox{for}\ x\in\mathcal{N}(\partial\Omega) $$
due to $f(x)\geq0$ for any $x\in\mathcal{N}(\partial\Omega)$.

Noting $v_\lambda(x)\gneqq 0$ on $\partial\mathcal{N}(\partial\Omega)$, by strong maximum principle, we have $v_\lambda(x)>0$ in $\mathcal{N}(\partial\Omega)$ for any $\lambda\in(0,\lambda_f)$. In conclusion, we have $v_\lambda(x)>0$ in $\Omega$ for any $\lambda\in(0,\lambda_f)$.
This completes the proof of the conclusion stated in {\it Step3}.

Finally, combining the statements of {\it Step1}, {\it Step2} and {\it Step3} together, we reach Lemma 3.3.
\vskip 0.1in
{\em {\bf Lemma 3.4}\ \ If $f(x)\in\mathcal{M}$ and $2\alpha+1<p<2^*$, then there exists a positive number $\lambda_0$ such that, for any $\lambda\in(0, \lambda_0)$, problem (\ref{eq3}) has at least one positive solution $u_\lambda$ with property $J_\lambda(u_\lambda)<0$.}
\vskip 0.1in
{\bf Proof:}\ \ Let $\rho_0$, $E_0$ and $\lambda_f$ be positive numbers determined in Lemma 3.3. Set $B_{\rho_0}=\{\ u\in H_0^1(\Omega):\ \ \|\nabla u\|_2<\rho_0\ \}$, and define
$$C_\lambda=\inf\limits_{u\in B_{\rho_0}}J_\lambda(u).$$
we can claim that $C_\lambda<0$. In fact, by the assumption $f(x)\in\mathcal{M}$, we know that problem (\ref{eq21}) has a solution $\varphi_\lambda$ which satisfies
$$\|\nabla\varphi_\lambda\|_2^2+b\|\nabla\varphi_\lambda\|_2^{2\alpha}=\lambda\int_\Omega f(x)\varphi_\lambda dx.$$
From this we can infer that
$$\|\nabla\varphi_\lambda\|_2\leq\frac{\sqrt{2}\|f\|_2}{\sqrt{\lambda_1(\Omega)}}\lambda.$$
Therefore, if we choose $\lambda_*=\frac{\sqrt{\lambda_1(\Omega)}\rho_0}{2\|f\|_2}$, then 
$$\|\nabla\varphi_\lambda\|_2\leq\frac{\rho_0}{\sqrt{2}}\ \ \ \mbox{for any}\ \ \lambda\in (0, \lambda_*).$$
This implies that $\varphi_\lambda\in B_{\rho_0}$ for any $\lambda\in (0, \lambda_*)$.  Noting $\alpha>0$, we have
$$J_\lambda(\varphi_\lambda)=-\frac{1}{2}\|\nabla\varphi_\lambda\|_2^2-b\frac{2\alpha+1}{2(\alpha+1)}\|\nabla\varphi_\lambda\|_2^{2(\alpha+1)}-\frac{1}{p+1}\|\varphi_\lambda\|_{p+1}^{p+1}<0.$$
By the definition of $C_\lambda$, we have
$$C_\lambda\leq J_\lambda(\varphi_\lambda)<0.$$
Let $\lambda_0=\min\{\lambda_f, \lambda_*\}$. For any fixed $\lambda\in (0, \lambda_0)$, if $\{u_{\lambda, n}\}$ is a minimizing sequence of $C_\lambda$, then we can claim that
$$\|\nabla u_{\lambda, n}\|_2\leq\rho_1$$
for some positive constant $\rho_1< \rho_0$. Otherwise, up to a subsequence, we may assume
$$\lim\limits_{n\rightarrow +\infty}\|\nabla u_{\lambda, n}\|_2=\rho_0.$$
Since $\lambda\in (0, \lambda_0)$, and
$$J_\lambda(u_{\lambda, n})\geq\frac{1}{4}\|\nabla u_{\lambda, n}\|^2_2-\frac{S(\Omega)}{p+1}\|\nabla u_{\lambda, n}\|^{p+1}_2-\frac{\lambda^2}{\lambda_1(\Omega)}\|f\|^2_2,$$
we have
$$0>C_\lambda=\lim\limits_{n\rightarrow +\infty}J_\lambda(u_{\lambda, n})\geq\frac{1}{4}\rho_0^2-\frac{S(\Omega)}{p+1}\rho_0^{p+1}-\frac{\lambda^2}{\lambda_1(\Omega)}\|f\|^2_2\geq E_0>0.$$
This is a contradiction.

By Ekeland's variational principle, we can find a sequence $\{v_{\lambda, n}\}$ such that
\begin{equation}\label{31}
\begin{array}{ll}
\lim\limits_{n\rightarrow +\infty}\|\nabla(u_{\lambda, n}-v_{\lambda, n})\|= 0,\\
\lim\limits_{n\rightarrow +\infty}J_\lambda(v_{\lambda, n})=C_\lambda,\\
\lim\limits_{n\rightarrow +\infty}J_\lambda^\prime(v_{\lambda, n})=0.
\end{array}
\end{equation}
Since $2\alpha+1<p<2^*$, a similar argument to that used in the proof of Theorem 2.1 implies that, up to a subsequence, $\{v_{\lambda, n}\}$ converges in $H_0^1(\Omega)$ to a function $v\in H_0^1(\Omega)$. Moreover, by a similar argument used in the proof of Lemma 3.3, we can prove that $v$ is a positive solution of  problem (\ref{eq3}). This completes the proof of Lemma 3.4.

To prove the nonexistence part of Theorem 3.1, we need the following result about semilinear problem
\vskip 0.1in
{\em {\bf Lemma 3.5}(\cite{DG, DY})\ \ If $1<p<2^*$, or $p>2^*$ and $\Omega$ is star-shaped, then, for any $f(x)\in\mathcal{M}$, there exists a positive number $\lambda_f$ such that the semilinear problem
\begin{equation}\label{eq324}
\left\{\begin{array}{ll}
-\Delta u=u^p+\lambda f & x\in\Omega\\
u\geq 0  & x\Omega\\
u=0 & x\in\partial\Omega
\end{array}
\right.
\end{equation}
has at least one solution for $\lambda\in (0, \lambda_f)$, and has no solution for $\lambda>\lambda_f$. Moreover, there exist a positive constant $C$ independent of $\lambda$ such that for any solution $u_\lambda$ of problem (\ref{eq324}) with respect to parameter $\lambda\in (0, \lambda_f)$, there holds
$$\|\nabla u_\lambda\|_2\leq C.$$ }

\vskip 0.1in
The nonexistence part of Theorem 3.1 is a special case of the following lemma.
\vskip 0.1in
{\em {\bf Lemma 3.6}\ \ If $2\alpha+1<p<2^*$, or $p>2^*$ and $\Omega$ is star-shaped, then, for any $f(x)\in\mathcal{M}$, there exists a positive number $\Lambda_f$ such that
problem (\ref{eq3}) has no positive solution for any $\lambda>\Lambda_f$.}
\vskip 0.1in
{\bf Proof:}\ \ If problem (\ref{eq3}) has a nonnegative solution $u_\lambda$ with respect to parameter $\lambda$, then we can see that  $v=\frac{u_\lambda}{(1+b\|\nabla u_\lambda\|^{2\alpha}_2)^{\frac{1}{p-1}}}$ is a solution of 

\begin{equation}\label{eq325}\begin{split}
\left\{\begin{array}{ll}
-\Delta v=v^p+\frac{\lambda}{(1+b\|\nabla u_\lambda\|^{2\alpha}_2)^{\frac{p}{p-1}}}f & \mbox{in}\ \ \Omega, \\
v\geq 0  &\mbox{in}\ \ \Omega,\\
v=0 & \mbox{on}\ \ \partial\Omega.
\end{array}
\right.
\end{split}\end{equation}
Therefore, by Lemma 3.5, we should have
\begin{equation}\label{eq326}
\lambda\leq\lambda_f(1+b\|\nabla u_\lambda\|^{2\alpha}_2)^{\frac{p}{p-1}}
\end{equation}
with $\lambda_f$ being the fixed number given in Lemma 3.5.

Furthermore, by the definition of $v$ and Lemma 3.5, we see that the following inequality hold for absolute positive constant $C$ given in Lemma 3.5.
$$\|\nabla u_\lambda\|_2^{p-1}\leq C(1+b\|\nabla u_\lambda\|^{2\alpha}_2)=bC\|\nabla u_\lambda\|^{2\alpha}_2+C.$$
Noting $p-1>2\alpha$, we can conclude from the above inequality that 
\begin{equation}\label{eq327}
\|\nabla u_\lambda\|_2\leq C
\end{equation}
for some positive constant $C$ independent of $\lambda$.

Substituting (\ref{eq327}) into (\ref{eq326}), we get
$$\lambda\leq\lambda_f(1+bC^{2\alpha})^{\frac{p}{p-1}}.$$
This implies that problem (\ref{eq3}) has no positive solution for $\lambda>\Lambda_f=\lambda_f(1+bC^{2\alpha})^{\frac{p}{p-1}}$. Therefore, the proof of Lemma 3.6 is completed.
\vskip 0.1in
{\bf Proof of Theorem 3.1:}\ \ If $f(x)\in\mathcal{M}$ and $2\alpha+1<p<2^*$, then it follows easily from Lemma 3.3 and Lemma 3.4 that there exists a positive number $\lambda_f$ such that problem (\ref{eq3}) has at least two positive solutions $u_\lambda$ and $v_\lambda$ with property $J_\lambda(u_\lambda)<0$ and $J_\lambda(v_\lambda)>0$ for any $\lambda\in (0, \lambda_f)$. The nonexistence part of Theorem 3.1 follows directly from Lemma 3.6. Therefore, we complete the proof of Theorem 3.1.
\section*{4. The case $p>2^*$}

\setcounter{section}{4}

\setcounter{equation}{0}

\noindent
In this section, we investigate the case $p>2^*$, and aim to proving the following theorem
\vskip 0.1in
{\em {\bf Theorem 4.1}\ \ If $p>2^*$ and $\Omega$ is starshaped, then for any $f(x)\in C^1(\Omega)\backslash\{0\}$ there are two positive number $\lambda_f$ and $\Lambda_f$ such that problem (\ref{eq3}) has at least one positive solution for any $\lambda\in (0, \lambda_f)$ if and only if $f(x)\in\mathcal{M}$, and has no positive solution for $\lambda>\Lambda_f$.}
\vskip 0.1in
Since $p>2^*$, we can not use variational method to get positive solution for problem (\ref{eq3}). At the same time, comparison principle may cease to validate for Kirchhoff type equations (see \cite{GFA}), we are also lack of sub-supersolution method for Kirchhoff type equation itself. Hence, some new ideas are needed for finding positive solutions of problem (\ref{eq3}) in this supercritical case. Here, we propose an iterative method based on the comparison principle of Laplace operator. The iterative sequence is no more monotone, but is still bounded. This is presented in the following lemma.
\vskip 0.1in
{\em {\bf Lemma 4.2}\ \ If $f(x)\in\mathcal{M}$, then there exists a positive number $\lambda_f$ such that problem (\ref{eq3}) has at least one positive solution for any $\lambda\in (0, \lambda_f)$.}
\vskip 0.1in
{\bf Proof:}\ \ Since $f(x)\in\mathcal{M}$, we can easily see that, for any $\lambda>0$, the following problem has a solution $\varphi_\lambda(x)$.
\begin{equation}\label{41}
\left\{
\begin{array}{ll}
-\Delta\varphi=\lambda f(x) & x\in\Omega,\\
\varphi\geq 0 & x\in\Omega,\\
\varphi=0 & x\in\partial\Omega.
\end{array}
\right.
\end{equation}
Let $\psi(x)$ be the solution of the following problem
\begin{equation}\label{42}
\left\{
\begin{array}{ll}
-\Delta\psi=1 & x\in\Omega,\\
\psi=0 & x\in\partial\Omega.
\end{array}
\right.
\end{equation}
Choosing $M_0>0$ so small that
$$M_0>M_0^p\max\limits_{x\in\Omega}\psi^p(x)+M_0^p\max\limits_{x\in\Omega}\|f(x)\|,$$
and setting $\psi_0(x)=M_0\psi(x)$, we can easily check that
\begin{equation}\label{43}
\left\{
\begin{array}{ll}
-\Delta\psi_0=M_0\geq\psi_0^p+\lambda f(x) & x\in\Omega\\
\psi_0=0 & x\in\partial\Omega
\end{array}
\right.
\end{equation}
for any $\lambda\in (0, M_0^p)$.

Taking (\ref{41}) and (\ref{43}) into account, we infer from the strong comparison principle for Laplace operator that
\begin{equation}\label{44}
\varphi_\lambda(x)<\psi_0(x)\ \ \mbox{for}\ \ x\in\Omega\ \ \mbox{and}\ \ \lambda\in (0, M_0^p).
\end{equation} 
Let $\lambda_f=M_0^p$. To obtain a solution of problem (\ref{eq3}) for any $\lambda\in (0, \lambda_f)$, we construct an approximation sequence $\{u_n(x)\}_{n=1}^\infty$ in the following way. 

Initially, we set $u_1(x)=\varphi_\lambda(x)$. Then, we get $u_{n+1}(x)$ from $u_n(x)$ by solving the following problem
\begin{equation}\label{45}
\left\{
\begin{array}{ll}
-(1+b\|\nabla u_{n+1}\|_2^{2\alpha})\Delta u_{n+1}=u_n^p+\lambda f(x) & x\in\Omega,\\
u_{n+1}=0 & x\in\partial\Omega.
\end{array}
\right.
\end{equation}
By induction method, we can see that
\begin{equation}\label{46}
0\leq u_n(x)\leq\psi_0(x)\ \ \mbox{for}\ \ x\in\Omega\ \ \mbox{and}\ \ n=1, 2, \cdots.
\end{equation}
Indeed, from (\ref{44}) , we firstly have
$$0\leq u_1(x)=\varphi_\lambda(x)\leq\psi_0(x).$$
If we inductively assume 
\begin{equation}\label{47}
0\leq u_k(x)\leq\psi_0(x).
\end{equation}
then what we should do is to proving 
\begin{equation}\label{48}
0\leq u_{k+1}(x)\leq\psi_0(x).
\end{equation}
Obviously, (\ref{48}) can be deduced from (\ref{47}) and the comparison principle of Laplace operator. In fact, on one hand, (\ref{41}), (\ref{45}) and (\ref{47}) imply that
\begin{equation}\label{49}
\left\{
\begin{array}{ll}
-(1+b\|\nabla u_{k+1}\|_2^{2\alpha})\Delta u_{k+1}=u_k^p+\lambda f(x)\geq\lambda f(x)=-\Delta\varphi_\lambda, & x\in\Omega,\\
u_{k+1}=\varphi_\lambda=0 & x\in\partial\Omega.
\end{array}
\right.
\end{equation}
Therefore, it follows from the comparison principle of Laplace operator that
\begin{equation}\label{410}
u_{k+1}(x)\geq\frac{\varphi_\lambda(x)}{1+b\|\nabla u_{k+1}\|_2^{2\alpha}}\geq 0.
\end{equation}
On the other hand, (\ref{43}), (\ref{45}) and (\ref{47}) imply
\begin{equation}\label{411}
\left\{
\begin{array}{ll}
-(1+b\|\nabla u_{k+1}\|_2^{2\alpha})\Delta u_{k+1}=u_k^p+\lambda f(x)\leq\psi_0^p+\lambda f(x)\leq-\Delta\psi_0, & x\in\Omega,\\
u_{k+1}=\psi_0=0 & x\in\partial\Omega.
\end{array}
\right.
\end{equation}
Hence, by the comparison principle of Laplace operator, we have
\begin{equation}\label{412}
u_{k+1}(x)\leq\frac{\psi_0(x)}{1+b\|\nabla u_{k+1}\|_2^{2\alpha}}\leq\psi_0(x).
\end{equation}
Combining (\ref{410}) and (\ref{412}) together, we get (\ref{48}). This concludes (\ref{46}) by induction method.

With (\ref{46}) established, we can deduce from (\ref{45}) and (\ref{43}) that
$$\|\nabla u_{n+1}\|_2\leq\|\nabla\psi_0\|_2.$$
From this and Schaulder's estimate, we have
$$\|u_{n+1}\|_{C^{2, \tau}(\Omega)}\leq C$$
for some positive constant $C$ and $\tau\in (0, 1)$ independent of $n$. Therefore, up to a subsequence, we may conclude that $\{u_n\}$ converges in $C^2(\Omega)$ to a function $u$ which is obviously a nonnegative solution of problem (\ref{eq3}). The positivity of $u$ follows from the strong comparison principle of Laplace operator. This completes the proof of Lemma 4.2.

The necessarity part of Theorem 4.1 includes in the following lemma

\vskip 0.1in
{\em {\bf Lemma 4.3}\ \ Assume that $p>2^*$, $f(x)\in C^1(\Omega)\backslash\{0\}$ and $\Omega$ is starshaped. If there exists a positive number $\lambda_f$ such that problem (\ref{eq3}) has positive solution for any $\lambda\in (0, \lambda_f)$, then $f(x)\in\mathcal{M}$.}
\vskip 0.1in
To prove Lemma 4.3, we need the following well known Pohozaev identity.
\vskip 0.1in
{\em {\bf Lemma 4.4 }(\cite{FLN}) Let $\Omega$ be a smooth bounded domain and suppose that $g: \Omega\times R\rightarrow R$ is a continuous map and that $\omega\in C^2(\overline{\Omega})$ satisfies
\begin{equation}\label{eq57}\begin{split}
\left\{\begin{array}{ll}
\Delta \omega+g(x,\omega(x))=0 & \mbox{in}\ \ \Omega, \\
\omega=0 & \mbox{on}\ \ \partial\Omega.
\end{array}
\right.
\end{split}\end{equation}
If $\nu(x)$ denotes the unit outward normal to $\partial\Omega$ at $x$, then $\omega$ satisfies
\begin{equation}
\int_{\partial\Omega}x\cdot\nu(x)|\nabla \omega|^2 dS 
=2N\int_\Omega G(x,\omega)dx+2\int_\Omega x\cdot \nabla_x G dx-(N-2)\int_\Omega g(x,\omega) \omega dx.
\end{equation}
where $G(x,\omega)=\int^\omega_0 g(x,t)dt$, and $\nabla_x G(x,\omega)$ is the gradient of $G(x,\omega)$ with respect to the variable $x$.}
\vskip 0.1in
{\bf Proof of Lemma 4.3:}\ \ No loss of generality, we may assume that $\Omega$ is star-shaped with respect to the origin $O$. That is $x\cdot\nu(x)\geq 0$ for any $x\in\partial\Omega$. Let $u_\lambda$ be positive solution of problem (\ref{eq3}) with respect to parameter $\lambda\in (0, \lambda_f)$. Setting
$$u_\lambda=\frac{\lambda}{1+b\|\nabla u_\lambda\|^{2\alpha}_2}v_\lambda, $$
we see that $v_\lambda$ satisfies

\begin{equation}\label{eq41}\begin{split}
\left\{\begin{array}{ll}
-\Delta v_\lambda=\frac{\lambda^{p-1}}{(1+b\|\nabla u_\lambda\|^{2\alpha}_2)^p}v^p_\lambda + f(x) & \mbox{in}\ \ \Omega, \\
v_\lambda=0 & \mbox{on}\ \ \partial\Omega.
\end{array}
\right.
\end{split}\end{equation}
Applying Lemma 4.4 to problem (\ref{eq41}), we have
$$\int_{\partial\Omega}x\cdot\nu(x) |\nabla v_\lambda|^2 dS=
\frac{-\eta\lambda^{p-1}}{(1+b\|\nabla u_\lambda\|^{2\alpha}_2)^p}\int_\Omega v^{p+1}_\lambda dx+2\int_\Omega x\cdot\nabla fv_\lambda dx+(2+N)\int_\Omega fv_\lambda dx$$
with $\eta=N-2-\frac{2N}{p+1}$. it worth mentioning here that $\eta>0$ due to $p>2^*$.

Since $\Omega$ is star-shaped with respect to $O$, we have
$$\int_{\partial\Omega}x\cdot\nu(x) |\nabla v_\lambda|^2 dS\geq 0.$$
Therefore
\begin{equation}\label{eq42}
\frac{\lambda^{p-1}}{(1+b\|\nabla u_\lambda\|^{2\alpha}_2)^p}\int_\Omega v^{p+1}_\lambda dx
\leq\frac{2}{\eta}\int_\Omega x\cdot\nabla fv_\lambda dx+\frac{N+2}{\eta}\int_\Omega fv_\lambda dx.
\end{equation}
By (\ref{eq41}), we have
\begin{align}\label{eq43}
\int_\Omega|\nabla v_\lambda|^2 dx=\frac{\lambda^{p-1}}{(1+b\|\nabla u_\lambda\|^{2\alpha}_2)^p}\int_\Omega v^{p+1}_\lambda dx+\int_\Omega fv_\lambda dx.
\end{align}
Combining (\ref{eq42}) and (\ref{eq43}) together, we get
$$\int_\Omega|\nabla v_\lambda|^2 dx\leq \frac{2}{\eta}\int_\Omega x\cdot\nabla fv_\lambda dx+(1+\frac{N+2}{\eta})\int_\Omega fv_\lambda dx.$$
This implies that
\begin{align}\label{eq44}
\int_\Omega|\nabla v_\lambda|^2 dx\leq C
\end{align}
for some positive constant $C$ independent of $\lambda$. 

Therefore, up to a subsequence, we have
$$v_\lambda(x)\rightharpoonup v(x) \ \ \mbox{weakly in}\ \ H_0^1(\Omega)\ \ \mbox{as}\ \ \lambda\rightarrow 0.$$
That is
\begin{align}\label{eq45}
\int_\Omega\nabla v_\lambda\cdot\nabla\varphi dx\rightarrow\int_\Omega\nabla v\cdot\nabla\varphi dx\ \ \mbox{as}\ \ \lambda\rightarrow 0
\end{align}
for any $\varphi(x)\in C_0^\infty(\Omega)$.

By (\ref{eq43})) and (\ref{eq44}), we have
\begin{align}\label{eq46}
\frac{\lambda^{p-1}}{(1+b\|\nabla u_\lambda\|^{2\alpha}_2)^p}\int_\Omega v^{p+1}_\lambda dx\leq C
\end{align}
for some positive constant $C$ independent of $\lambda$. Consequently, for any $\varphi(x)\in C_0^\infty(\Omega)$, we have
\begin{align}\label{eq47}
\frac{\lambda^{p-1}}{(1+b\|\nabla u_\lambda\|^{2\alpha}_2)^p}\int_\Omega v^p_\lambda\varphi dx\rightarrow 0\ \ \mbox{as}\ \ \lambda\rightarrow 0.
\end{align}
By (\ref{eq41}), for any $\varphi(x)\in C_0^\infty(\Omega)$, we have
\begin{align}\label{eq48}
\int_\Omega\nabla v_\lambda\nabla\varphi dx=\frac{\lambda^{p-1}}{(1+b\|\nabla u_\lambda\|^{2\alpha}_2)^p}\int_\Omega v^p_\lambda\varphi dx+\int_\Omega f\varphi dx.
\end{align}
Sendding $\lambda$ to $0$ in (\ref{eq48}), and taking (\ref{eq45}) and (\ref{eq47}) account, we get
$$\int_\Omega\nabla v\nabla\varphi dx=\int_\Omega f\varphi dx,\ \ \mbox{for any}\ \ \varphi\in C^\infty_0(\Omega).$$
This and the regularity theory of elliptic equations imply that $v$ is a solution of problem (\ref{eq5}). Moreover, $v(x)\geq 0$ due to $v_\lambda(x)\geq 0$. Therefore, $f(x)\in\mathcal{M}$. This completes the proof of Lemma 4.3.
\vskip 0.1in

{\bf Proof of Theorem 4.1}\ \  Combining Lemma 3.6, Lemma 4.2 and Lemma 4.3 together, we reach the conclusion of Theorem 4.1.

\end{document}